\documentclass[]{ifacconf}
\usepackage{multirow}

\usepackage{epsfig} 
\usepackage{amsmath} 
\usepackage{amssymb}  
\usepackage{url}

\usepackage{enumerate}

\begin{document}

\begin{frontmatter}

\title{Improving Fast Dual Ascent for MPC - Part II: The Embedded Case\thanksref{footnoteinfo}}

\thanks[footnoteinfo]{During the preparation of this paper, the author was a
member of the LCCC Linnaeus Center at Lund University. 
Financial support from the Swedish Research Council for the author's
Postdoctoral studies at Stanford University is gratefully
acknowledged. Further, Eric Chu is gratefully acknowledged for
constructive feedback and Alexander Domahidi is gratefully
acknowledged for suggesting the AFTI-16 control problem as a benchmark
and providing FORCES code for the same.}

\author[Stanford]{Pontus Giselsson} 

\address[Stanford]{Electrical Engineering, Stanford University\\
   (e-mail: pontusg@stanford.edu).}                                              




\begin{abstract}

Recently, several authors have suggested the use of first order
methods, such as fast dual ascent and the alternating direction method
of multipliers, for embedded model predictive control. The main reason is that
they can be implemented using simple arithmetic operations only.
However, a known limitation of gradient-based methods is that they are
sensitive to ill-conditioning of the problem data. In this paper, we
present a fast dual gradient method for which the sensitivity to
ill-conditioning is greatly reduced. This is achieved by approximating
the negative dual function with
a quadratic upper bound with different curvature in different
directions in the algorithm, as
opposed to having the same curvature in all directions as in standard
fast gradient methods. The main contribution of this paper is a
characterization of the set of matrices that can be used to form such
a quadratic upper bound to the negative dual function. We also
describe how to choose
a matrix from this set to get an improved approximation of the dual
function, especially if it is ill-conditioned, compared to the approximation used in standard
fast dual gradient methods. This can give a significantly improved
performance as illustrated by a numerical
evaluation on an ill-conditioned AFTI-16 aircraft model.

\end{abstract}

\end{frontmatter}




\section{Introduction}

Several authors including
\cite{ODonoghueSplitting,JerezMHz,RichterMMOR2013,PatrinosTAC2013}
have recently proposed first order optimization
methods as appropriate for embedded model predictive control. In
\cite{ODonoghueSplitting,JerezMHz},
the alternating direction method of multipliers (ADMM, see
\cite{BoydDistributed}) were used and
high computational speeds were reported when implemented on embedded
hardware. In \cite{RichterMMOR2013,PatrinosTAC2013}, the optimal
control problems arising in model predictive were solved using different
formulations of fast dual gradient methods. In \cite{RichterMMOR2013},
the equality constraints, i.e. the dynamic constraints, are dualized and
a diagonal cost and box constraints are assumed. The resulting dual
problem is solved using a fast gradient method. In
\cite{PatrinosTAC2013}, the same splitting as in
\cite{ODonoghueSplitting,JerezMHz} is used, but a fast gradient method
is used to solve the resulting problem as opposed to ADMM in
\cite{ODonoghueSplitting,JerezMHz}. In this paper, we will show how to
improve and generalize the fast dual gradient methods presented in
\cite{RichterMMOR2013,PatrinosTAC2013}.

Fast gradient methods as used in \cite{RichterMMOR2013,PatrinosTAC2013} have
been around since the early 80's
when the seminal paper \cite{Nesterov1983} was published. However,
fast gradient methods did not render much attention before the mid
00's, after which an increasing interest has emerged. Several
extensions and generalizations of the fast gradient method have been
proposed, e.g. in \cite{NesterovLectures,Nesterov2005}.
In \cite{BecTab_FISTA:2009}, the method was generalized to
allow for minimization of composite objective functions. Further, a unified
framework for fast gradient methods and their generalizations were
presented in \cite{Tseng_acc:2008}.
To use fast gradient methods for composite minimization, one
objective term should be convex and differentiable with a Lipschitz continuous
gradient, while the other should be proper, closed, and convex. The
former condition is equivalent to the existence of a quadratic
upper bound to the function, with
the same curvature in all directions. The curvature is specified by
the Lipschitz constant to the gradient. In fast gradient methods, the
quadratic upper bound 
serves as an approximation of the function to be minimized, since the
bound is minimized in every iteration of the algorithm.
If the
quadratic upper bound does not well approximate the function to be
minimized, slow convergence properties are expected. By instead
allowing for a quadratic upper bound with different curvature in
different directions, as in generalized fast gradient methods
\cite{WangmengTIP2011}, the bound can closer approximate the function to be
minimized. For an appropriate choice of non-uniform quadratic upper bound, this
can significantly improve the performance of the algorithm.

In \cite[Theorem 1]{Nesterov2005}, a Lipschitz constant to the gradient of
the dual function to strongly convex problems is presented. This
result quantifies
the curvature of a uniform quadratic upper bound to the negative dual
function. This result was improved in \cite[Theorem
7]{RichterMMOR2013} when 
the primal cost is restricted to being quadratic. Using these
quadratic upper bounds, with the same curvature in all directions, as
dual function approximation in a fast dual gradient method, may result in slow
convergence rates. Especially for ill-conditioned problems
where the upper bound does not well approximate the negative dual
function. In this paper, the main result is a new characterization of 
the set of matrices that can be used to describe quadratic upper
bounds to the negative dual function.
This result generalizes and improves previous
results in \cite{Nesterov2005,RichterMMOR2013}. We also show how to
appropriately choose a matrix from this set to get a quadratic upper
bound that well approximates the negative dual function. Since in the
proposed method, the dual function approximation is better that in
standard fast dual gradient methods used in
\cite{RichterMMOR2013,PatrinosTAC2013}, better convergence rate
properties are expected.

In model predictive control, much offline computational effort can be
devoted to improve the online execution time of the solver. This is
done, e.g., in explicit MPC, see \cite{Bemporad}, where the explicit
parametric solution is 
computed beforehand, and found through a look-up table online. In this
paper, the offline computational effort is devoted to choose a matrix
that describes the quadratic upper bound to the negative dual
function.
The computed matrix is the same in all samples in the controller
and can therefore be computed offline.
The algorithm is evaluated on a pitch control problem in an AFTI-16 aircraft
that has previously been studied in \cite{Kapasouris,Bemporad97}. This
is a challenging problem for first order methods since it is very
ill-conditioned. The numerical evaluation shows that the method
presented in this paper outperforms other first-order methods presented in
\cite{ODonoghueSplitting,JerezMHz,RichterMMOR2013,PatrinosTAC2013}
with one to three orders of magnitude. Also, the numerical evaluation
shows that a C implementation of our algorithm outperform FORCES,
\cite{FORCES}, which is a C code-generator for MPC problems using
a tailored interior point method, and the general commercial QP-solver
MOSEK.

This paper extends the conference publication \cite{gis2014IFACopt}, and is the second
of a series of two papers on improving duality-based optimization in
MPC, with \cite{gis2014AutPart1} being the first.

\section{Preliminaries and Notation}

\subsection{Notation}

We denote by $\mathbb{R}$, $\mathbb{R}^n$, $\mathbb{R}^{m\times n}$,
the sets of real numbers, vectors, and matrices.
$\mathbb{S}^n\subseteq\mathbb{R}^{n\times n}$ is the set of symmetric matrices, and 
$\mathbb{S}_{++}^n\subseteq\mathbb{S}^n$,
$[\mathbb{S}_{+}^n]\subseteq\mathbb{S}^n$, are the sets of positive [semi] 
definite matrices. Further, $L\succeq M$ and $L\succ M$
where $L,M\in\mathbb{S}^n$ denotes
$L-M\in\mathbb{S}_{+}^n$ and $L-M\in\mathbb{S}_{++}^n$ respectively.
We also use notation
$\langle x,y\rangle=x^Ty$, $\langle
x,y\rangle_H=x^THy$, $\|x\|_2=\sqrt{x^Tx}$, and $\|x\|_H =
\sqrt{x^THx}$. Finally, $I_{\mathcal{X}}$ denotes the indicator
function for the set $\mathcal{X}$, i.e.
$I_{\mathcal{X}}(x)\triangleq\left\{\begin{smallmatrix}0,~&x\in\mathcal{X}\\
\infty, &{\rm{else~}}\end{smallmatrix}\right.$.

\subsection{Preliminaries}

In this section, we introduce generalizations of well used
concepts. We generalize the notion of strong convexity as well as
the notion of Lipschitz continuity of the gradient of convex
functions. We also define conjugate
functions and state a known result on dual properties of a
function and its conjugate.

For differentiable and convex functions
$f~:~\mathbb{R}^n\to\mathbb{R}$ that have a Lipschitz continuous
gradient with constant $L$, we have that
\begin{equation}
\|\nabla f(x_1)-\nabla f(x_2)\|_2\leq L\|x_1-x_2\|_2
\label{eq:gradLipschitz}
\end{equation}
holds for all $x_1,x_2\in\mathbb{R}^n$.
This is equivalent to that
\begin{equation}
f(x_1)\leq f(x_2)+\langle \nabla
f(x_2),x_1-x_2\rangle+\frac{L}{2}\|x_1-x_2\|_2^2
\label{eq:standardQuadBound}
\end{equation}
holds for all
$x_1,x_2\in\mathbb{R}^n$ \cite[Theorem 2.1.5]{NesterovLectures}.
In this paper, we allow for a generalized version of the quadratic
upper bound \eqref{eq:standardQuadBound} to $f$, namely that
\begin{equation}
f(x_1)\leq f(x_2)+\langle \nabla
f(x_2),x_1-x_2\rangle+\frac{1}{2}\|x_1-x_2\|_{\mathbf{L}}^2
\label{eq:quadBound}
\end{equation}
holds for all $x_1,x_2\in\mathbb{R}^n$ where
$\mathbf{L}\in\mathbb{S}_+^{n}$. The
bound \eqref{eq:standardQuadBound} is obtained by setting
$\mathbf{L}=LI$ in \eqref{eq:quadBound}. 
\begin{rem}
For concave functions $f$, i.e. where $-f$ is convex, the Lipschitz
condition \eqref{eq:gradLipschitz} is equivalent to that the following
quadratic lower bound
\begin{equation}
f(x_1)\geq f(x_2)+\langle \nabla
f(x_2),x_1-x_2\rangle-\frac{L}{2}\|x_1-x_2\|_2^2
\label{eq:standardQuadBoundConcave}
\end{equation}
holds for all $x_1,x_2\in\mathbb{R}^n$. The generalized counterpart
naturally becomes that 
\begin{equation}
f(x_1)\geq f(x_2)+\langle \nabla
f(x_2),x_1-x_2\rangle-\frac{1}{2}\|x_1-x_2\|_{\mathbf{L}}^2
\label{eq:quadBoundConcave}
\end{equation}
holds for all $x_1,x_2\in\mathbb{R}^n$.
\label{rem:quadBoundConcave}
\end{rem}
Next, we state a
Lemma on equivalent characterizations of the condition \eqref{eq:quadBound}.
\begin{lem}
Assume that $f~:~\mathbb{R}^n\to\mathbb{R}$
is convex and differentiable. The condition that 
\begin{equation}
f(x_1)\leq f(x_2)+\langle
\nabla f(x_2),x_1-x_2\rangle+\frac{1}{2}\|x_1-x_2\|_{\mathbf{L}}^2
\label{eq:quadBoundLem}
\end{equation}
holds for some
$\mathbf{L}\in\mathbb{S}_{+}^n$ and all $x_1,x_2\in\mathbb{R}^n$
is equivalent to that 
\begin{equation}
\langle\nabla f(x_1)-\nabla f(x_2),x_1-x_2\rangle \leq \|x_1-x_2\|_{\mathbf{L}}^2.
\label{eq:quadBoundCondition}
\end{equation}
holds for all $x_1,x_2\in\mathbb{R}^n$.
\label{lem:quadUBequiv}
\end{lem}
\begin{pf}
To show the equivalence, we introduce the function
$g(x):=\frac{1}{2}x^T\mathbf{L}x-f(x)$. According to \cite[Theorem
2.1.3]{NesterovLectures} and since $g$ is differentiable,
$g~:~\mathbb{R}^n\to\mathbb{R}$ is convex if and only if $\nabla g$ is
monotone. The function $g$ is convex if and only if 
\begin{align*}
&g(x_1)\geq g(x_2)+\langle\nabla g(x_2),x_1-x_2\rangle=\\
&=\frac{1}{2}x_2^T\mathbf{L}x_2-f(x_2)+\langle \mathbf{L}x_2-\nabla
f(x_2),x_1-x_2\rangle\\
&=-f(x_2)-\langle \nabla
f(x_2),x_1-x_2\rangle-\tfrac{1}{2}\|x_1-x_2\|_\mathbf{L}^2+\tfrac{1}{2}x_1^T\mathbf{L}x_1.
\end{align*}
Noting that $g(x_1) = \frac{1}{2}x_1^T\mathbf{L}x_1-f(x_1)$ gives the
negated version of \eqref{eq:quadBoundLem}.

Monotonicity of $\nabla g$ is equivalent to 
\begin{align*}
0&\leq \langle \nabla g(x_1)-\nabla g(x_2),x_1-x_2\rangle\\
&=\langle \mathbf{L}x_1-\nabla f(x_1)-\mathbf{L}x_2+\nabla f(x_2),x_1-x_2\rangle\\
&=\|x_1-x_2\|_{\mathbf{L}}^2-\langle \nabla f(x_1)-\nabla f(x_2),x_1-x_2\rangle.
\end{align*}
Rearranging the terms gives \eqref{eq:quadBoundCondition}. This
concludes the proof.
\end{pf}

Next, we state the corresponding result for concave functions.
\begin{cor}
Assume that $f~:~\mathbb{R}^n\to\mathbb{R}$
is concave and differentiable. The condition that 
\begin{equation}
f(x_1)\geq f(x_2)+\langle
\nabla f(x_2),x_1-x_2\rangle-\frac{1}{2}\|x_1-x_2\|_{\mathbf{L}}^2
\label{eq:quadBoundCorConcave}
\end{equation}
holds for some
$\mathbf{L}\in\mathbb{S}_{+}^n$ and all $x_1,x_2\in\mathbb{R}^n$
is equivalent to that 
\begin{equation}
\langle\nabla f(x_1)-\nabla f(x_2),x_2-x_1\rangle \leq \|x_1-x_2\|_{\mathbf{L}}^2.
\label{eq:quadBoundConditionConcave}
\end{equation}
holds for all $x_1,x_2\in\mathbb{R}^n$.
\label{cor:quadUBequivConcave}
\end{cor}
\begin{pf}
The proof follows directly from $-f$ being convex and applying
Lemma~\ref{lem:quadUBequiv}. 
\end{pf}

The standard definition of a differentiable and strongly convex 
function $f~:~\mathbb{R}^n\to\mathbb{R}$ is that it satisfies
\begin{equation}
f(x_1)\geq f(x_2)+\langle \nabla
f(x_2),x_1-x_2\rangle+\frac{\sigma}{2}\|x_1-x_2\|_2^2 
\label{eq:strConvStandardDef}
\end{equation}
for any $x_1,x_2\in\mathbb{R}^n$, where the modulus $\sigma\in\mathbb{R}_{++}$
describes a lower bound on the curvature of the function. In this
paper, the definition \eqref{eq:strConvStandardDef} is generalized to
allow for a quadratic lower bound with different curvature in different
directions.
\begin{defn}
A differentiable function $f~:~\mathbb{R}^n\to\mathbb{R}$
is \emph{strongly convex with matrix H} if and only if
\begin{equation*}
f(x_1)\geq f(x_2)+\langle \nabla
f(x_2),x_1-x_2\rangle+\frac{1}{2}\|x_1-x_2\|_H^2 
\end{equation*}
holds for all $x_1,x_2\in\mathbb{R}^n$, where
$H\in\mathbb{S}_{++}^{n}$.
\label{def:strConv}
\end{defn}
\begin{rem}
The traditional definition of strong convexity
\eqref{eq:strConvStandardDef} is obtained from
Definition~\ref{def:strConv} by setting $H=\sigma I$.
\end{rem}
\begin{lem}
Assume that $f~:~\mathbb{R}^n\to\mathbb{R}$
is differentiable and
strongly convex with matrix $H$. The condition that
\begin{equation}
f(x_1)\geq f(x_2)+\langle
\nabla f(x_2),x_1-x_2\rangle+\frac{1}{2}\|x_1-x_2\|_{H}^2
\label{eq:strConv}
\end{equation}
holds for all $x_1,x_2\in\mathbb{R}^n$
is equivalent to that
\begin{equation}
\langle\nabla f(x_1)-\nabla f(x_2),x_1-x_2\rangle \geq \|x_1-x_2\|_{H}^2
\label{eq:strConvCondition}
\end{equation}
holds for all $x_1,x_2\in\mathbb{R}^n$.
\label{lem:strConvCondition}
\end{lem}
\begin{pf}
To show the equivalence, we introduce the function
$g(x):=f(x)-\frac{1}{2}x^THx$ and proceed similarly to in the proof
of Lemma \eqref{lem:quadUBequiv}. According to \cite[Theorem
2.1.3]{NesterovLectures} and since $g$ is differentiable,
$g~:~\mathbb{R}^n\to\mathbb{R}$ is convex if and only if $\nabla g$ is
monotone. The function $g$ is convex if and only if 
\begin{align*}
&g(x_1)\geq g(x_2)+\langle\nabla g(x_2),x_1-x_2\rangle=\\
&=f(x_2)-\frac{1}{2}x_2^THx_2+\langle \nabla
f(x_2)-Hx_2,x_1-x_2\rangle\\
&=f(x_2)+\langle \nabla
f(x_2),x_1-x_2\rangle+\tfrac{1}{2}\|x_1-x_2\|_H^2-\tfrac{1}{2}x_1^THx_1.
\end{align*}
Noting that $g(x_1) = f(x_1)-\frac{1}{2}x_1^THx_1$ gives \eqref{eq:strConv}.

Monotonicity of $\nabla g$ is equivalent to 
\begin{align*}
0&\leq \langle \nabla g(x_1)-\nabla g(x_2),x_1-x_2\rangle\\
&=\langle \nabla f(x_1)-Hx_1-\nabla f(x_2)+\mathbf{L}x_2,x_1-x_2\rangle\\
&=\langle \nabla f(x_1)-\nabla f(x_2),x_1-x_2\rangle-\|x_1-x_2\|_{H}^2.
\end{align*}
Rearranging the terms gives \eqref{eq:strConvCondition}. This
concludes the proof.
\end{pf}

The condition \eqref{eq:strConv} is a quadratic lower bound on the function value,
while the condition \eqref{eq:quadBound} is a quadratic upper bound on the function
value. These two properties are linked through the conjugate function
\begin{align*}
f^{\star}(y) &\triangleq \sup_x \left\{y^Tx-f(x)\right\}.
\end{align*}
More precisely, we have the following result.
\begin{prop}
Assume that $f~:~\mathbb{R}^n\to\mathbb{R}\cup\{\infty\}$ is closed,
proper, and strongly convex with
modulus $\sigma$ on the relative interior of its domain. Then the
conjugate function
$f^{\star}$ is convex and differentiable, and
$\nabla f^{\star}(y) = x^\star(y)$, where $x^\star(y)=\arg\max_x 
\left\{y^Tx-f(x)\right\}$. Further, $\nabla f^{\star}$ is Lipschitz continuous with 
constant $L = \frac{1}{\sigma}$.
\label{prp:strConvSmooth}
\end{prop}
A straight-forward generalization is given by the chain-rule and was
proven in \cite[Theorem 1]{Nesterov2005} (which also proves the less
general Proposition~\ref{prp:strConvSmooth}).
\begin{cor}
Assume that $f~:~\mathbb{R}^n\to\mathbb{R}\cup\{\infty\}$ is closed,
proper, and strongly convex with
modulus $\sigma$ on the relative interior of its domain. Further, define
$g^\star(y) \triangleq f^\star(Ay)$. Then
$g^\star$ is convex and differentiable, and $\nabla g^\star(y) =
A^Tx^\star(Ay)$, where $x^\star(Ay)=\arg\max_x 
\left\{(Ay)^Tx-f(x)\right\}$. Further, $\nabla g^{\star}$
is Lipschitz continuous with 
constant $L = \frac{\|A\|_2^2}{\sigma}$.
\label{cor:strConvSmooth}
\end{cor}
For the case when $f(x) = \frac{1}{2}x^THx+g^Tx$, i.e. $f$ is a
quadratic, a tighter Lipschitz constant to $\nabla g^\star(y)=\nabla f^\star(Ay)$ was provided
in \cite[Theorem 7]{RichterMMOR2013}, namely
$L=\|AH^{-1}A^T\|_2$.

\section{Problem formulation}

We consider optimization problems of the form
\begin{equation}
\begin{tabular}[t]{ll}
minimize & $f(x)+h(x)+g(Bx)$\\
subject to & $Ax=b$
\end{tabular}
\label{eq:primProb}
\end{equation}
where $x\in\mathbb{R}^n$, $A\in\mathbb{R}^{m\times n}$,
$B\in\mathbb{R}^{p\times n}$,
$b\in\mathbb{R}^m$. We assume that the following assumption holds
throughout the paper:
\begin{assum}~
\begin{enumerate}[(a)]
\item The function $f~:~\mathbb{R}^n\to\mathbb{R}$ is differentiable
and strongly convex with matrix 
$H$.
\item The extended valued functions
$h~:~\mathbb{R}^n\to\mathbb{R}\cup\{\infty\}$ and
$g~:~\mathbb{R}^n\to\mathbb{R}\cup\{\infty\}$, are proper, closed, and
convex.
\item $A\in\mathbb{R}^{m\times 
  n}$ has full row rank.
\end{enumerate}
\label{ass:probAss}
\end{assum}
\begin{rem}
Examples of
functions that satisfy Assumption~\ref{ass:probAss}(a) and
\ref{ass:probAss}(b) are
$f(x)=\frac{1}{2}x^THx+g^Tx$ with $H\in\mathbb{S}_{++}^n$ for
Assumption~\ref{ass:probAss}(a), and
$g=I_{\mathcal{X}}$, $g=\|\cdot\|_1$, $g=I_{\mathcal{X}}^\star$, or $g=0$ for
Assumption~\ref{ass:probAss}(b). If Assumption~\ref{ass:probAss}(c) is
not satisfied, redundant equality constraints can, without affecting
the solution of \eqref{eq:primProb}, be removed to satisfy the assumption.
\end{rem}

The optimization problem \eqref{eq:primProb} can equivalently be
written as
\begin{equation}
\begin{tabular}[t]{ll}
minimize & $f(x)+h(x)+g(y)$\\
subject to & $Ax=b$\\
& $Bx=y$
\end{tabular}
\label{eq:primProbY}
\end{equation}
We introduce dual
variables $\lambda\in\mathbb{R}^m$ for the equality constraints $Ax=b$
and dual variables $\mu\in\mathbb{R}^p$ for the equality constraints
$Bx=y$. This gives the following Lagrange dual problem
\begin{align}
&\nonumber \displaystyle\sup_{\lambda,\mu}\inf_{x,y}\left\{
\displaystyle f(x)+h(x)+\lambda^T(Ax-b)+g(y)+\mu^T(Bx-y)\right\}\\
&\nonumber = \displaystyle\sup_{\lambda,\mu}\Big[-\sup_{x}\left\{
\displaystyle (-A^T\lambda-B^T\mu)^Tx-f(x)-h(x)\right\}\\
&\nonumber \qquad\qquad\qquad\qquad\qquad\qquad-b^T\lambda-\sup_{y} \left\{\mu^Ty-g(y)\right\}\Big]\\
&\displaystyle =\sup_{\lambda,\mu}\left\{-F^\star(-A^T\lambda-B^T\mu)-b^T\lambda-g^\star(\mu)\right\}&
\label{eq:dualProb}
\end{align}
where $F^\star$ is the conjugate function to $F:=f+h$ and $g^\star$ is
the conjugate function to $g$.
For ease of exposition, we introduce $\nu =
(\lambda,\mu)\in\mathbb{R}^{m+p}$, $C =
[A^T~B^T]^T\in\mathbb{R}^{(m+p)\times n}$, and $c =
(b,0)\in\mathbb{R}^{m+p}$ and the following function
\begin{equation}
d(\nu) := -F^\star(-C^T\nu)-c^T\nu=-F^\star(-A^T\lambda-B^T\mu)-b^T\lambda.
\label{eq:dFcnDef}
\end{equation}
This implies that the dual problem \eqref{eq:dualProb} can be written
as
\begin{equation}
{\hbox{maximize }} d(\nu)-g^\star([0~I]\nu).
\label{eq:dualProb_dg}
\end{equation}
To evaluate the function $d$, an optimization
problem is solved. The minimand to this problem is denoted by
\begin{align}
\label{eq:innerMinimizer} x^\star(\nu) &:= \arg\min_{x}\left\{F(x)+\nu^TCx\right\}  \\
\nonumber &= \arg\min_{x}\left\{f(x)+h(x)+\lambda^TAx+\mu^TBx\right\}.
\end{align}
From Corollary~\ref{cor:strConvSmooth} we get that the function $d$ is
concave and differentiable
with gradient
\begin{equation*}
\nabla d(\nu) =Cx^\star(\nu)-c
\end{equation*}
and that $\nabla d$ is Lipschitz continuous with constant
$L=\|C\|_2^2/\lambda_{\min}(H)$, i.e., that
\begin{equation}
\|\nabla d(\nu_1)-\nabla d(\nu_2)\|_2\leq L\|\nu_1-\nu_2\|_2
\label{eq:dLipschitz}
\end{equation}
holds for all $\nu_1,\nu_2\in\mathbb{R}^{m+p}$. As stated in
Remark~\ref{rem:quadBoundConcave}, \eqref{eq:dLipschitz} is
equivalent to that the following 
quadratic lower bound to
the concave function $d$ holds for all
$\nu_1,\nu_2\in\mathbb{R}^{m+p}$
\begin{equation*}
d(\nu_1)\geq d(\nu_2)+\langle
\nabla d(\nu_2),\nu_1-\nu_2\rangle-\frac{L}{2}\|\nu_1-\nu_2\|_2^2.
\end{equation*}
In the following section we will show that the function $d$
satisfies the following tighter condition
\begin{equation}
d(\nu_1)\geq d(\nu_2)+\langle
\nabla d(\nu_2),\nu_1-\nu_2\rangle-\frac{1}{2}\|\nu_1-\nu_2\|_{\mathbf{L}}^2
\label{eq:dualQuadBoundIntr}
\end{equation}
for all $\nu_1,\nu_2\in\mathbb{R}^{m+p}$ and
$\mathbf{L}\in\mathbb{S}_+^{m+p}$ that satisfies $\mathbf{L} \succeq
CH^{-1}C^T$. 

\section{Dual function properties}

To show that the function $d$ as defined in \eqref{eq:dFcnDef} satisfies
\eqref{eq:dualQuadBoundIntr}, we need the following lemma.
\begin{lem}
Suppose that Assumption~\ref{ass:probAss} holds. Then
\begin{equation*}
\|x^\star(\nu_1)-x^\star(\nu_2)\|_H\leq \|\nu_1-\nu_2\|_{CH^{-1}C^T}
\end{equation*}
holds for all $\nu_1,\nu_2\in\mathbb{R}^{m+p}$, where $x^{\star}(\nu)$ is
defined in \eqref{eq:innerMinimizer}.
\label{lem:helpResult}
\end{lem}
\begin{pf}
We first show that
\begin{multline}
\label{eq:helpResult} \langle\nabla f(x^\star(\nu_1))-\nabla
f(x^\star(\nu_2)),x^\star(\nu_1)-x^\star(\nu_2)\rangle\leq
\\
\langle C^T(\nu_1-\nu_2),x^\star(\nu_2)-x^\star(\nu_1)\rangle.
\end{multline}
First order optimality conditions of \eqref{eq:innerMinimizer}
with $\nu_1$ and $\nu_2$ respectively are
\begin{align}
\label{eq:firstOrdOpt1} 0&\in\nabla f(x^\star(\nu_1))+\partial
h(x^\star(\nu_1))+C^T\nu_1,\\
\label{eq:firstOrdOpt2} 0&\in\nabla f(x^\star(\nu_2))+\partial
h(x^\star(\nu_2))+C^T\nu_2.
\end{align}
We denote by $\xi(x^\star(\nu_1))\in\partial h(x^\star(\nu_1))$
and $\xi(x^\star(\nu_2))\in\partial h(x^\star(\nu_2))$ the
sub-gradients that give equalities in \eqref{eq:firstOrdOpt1} and
\eqref{eq:firstOrdOpt1} respectively. This gives
\begin{align}
\label{eq:firstOrdOptEq1} 0&=\nabla f(x^\star(\nu_1))+\xi(x^\star(\nu_1))+C^T\nu_1,\\
\label{eq:firstOrdOptEq2}0&=\nabla f(x^\star(\nu_2))+\xi(x^\star(\nu_2))+C^T\nu_2.
\end{align}
Taking the scalar product of \eqref{eq:firstOrdOptEq1} with 
$x^\star(\nu_2)-x^\star(\nu_1)$ and the scalar product of
\eqref{eq:firstOrdOptEq2} with 
$x^\star(\nu_1)-x^\star(\nu_2)$, and summing gives
\begin{align*}
\langle\nabla f(x^\star(\nu_1))-\nabla
f(x^\star(\nu_2)),x^\star(\nu_1)-x^\star(\nu_2)\rangle+\qquad&\\
+\langle
C^T(\nu_1-\nu_2),x^\star(\nu_1)-x^\star(\nu_2)\rangle&=\\
=\langle\xi(x^\star(\nu_1))-\xi(x^\star(\nu_2)),x^\star(\nu_2)-x^\star(\nu_1)\rangle&\leq
 0
\end{align*}
where the inequality holds since sub-differentials of proper,
closed, and convex functions are (maximal) monotone
mappings, see \cite[~\S 24]{Rockafellar}.
This implies that \eqref{eq:helpResult} holds.

Further
\begin{align*}
  \|&x^\star(\nu_1)-x^\star(\nu_2)\|_H^2\leq\\
&\leq \langle\nabla f(x^\star(\nu_1))-\nabla
f(x^\star(\nu_2)),x^\star(\nu_1)-x^\star(\nu_2)\rangle\\
&\leq \langle
C^T(\nu_1-\nu_2),x^\star(\nu_2)-x^\star(\nu_1)\rangle
\\
&=\langle
H^{-1/2}C^T(\nu_1-\nu_2),H^{1/2}(x^\star(\nu_2)-x^\star(\nu_1))\rangle\\
&\leq \|H^{-1/2}C^T(\nu_1-\nu_2)\|_2\|x^\star(\nu_2)-x^\star(\nu_1)\|_H
\end{align*}
where the first inequality comes from
Lemma~\ref{lem:strConvCondition}, the second from \eqref{eq:helpResult}, 
and the final inequality is due to Cauchy Schwarz inequality.
This implies that
\begin{equation*}
\|x^\star(\nu_1)-x^\star(\nu_2)\|_H\leq \|\nu_1-\nu_2\|_{CH^{-1}C^T}
\end{equation*}
which concludes the proof.
\end{pf}

Now we are ready to state the main theorem of this section.
\begin{thm}
The function $d$ defined in \eqref{eq:dFcnDef} is concave, differentiable and
satisfies 
\begin{equation}
d(\nu_1)\geq d(\nu_2)+\langle
\nabla d(\nu_2),\nu_1-\nu_2\rangle-\frac{1}{2}\|\nu_1-\nu_2\|_{\mathbf{L}}^2
\label{eq:dualQuadBound}
\end{equation}
for every $\nu_1,\nu_2\in\mathbb{R}^{m+p}$ and
$\mathbf{L}\in\mathbb{S}_{+}^{m+p}$ that satisfies $\mathbf{L}\succeq
CH^{-1}C^T$.
\label{thm:dualQuadBound}
\end{thm}
\begin{pf}
Concavity and differentiability is deduced from Danskin's Theorem, see
\cite[Proposition~B.25]{Bertsekas99}.

To show \eqref{eq:dualQuadBound}, we have for any
$\nu_1,\nu_2\in\mathbb{R}^{m+p}$ that
\begin{align*}
\langle\nabla d(\nu_1)-&\nabla
d(\nu_2),\nu_2-\nu_1\rangle=\\
&= \langle
Cx^\star(\nu_1)-c-Cx^\star(\nu_2)+c,\nu_2-\nu_1\rangle\\
&=\langle
x^\star(\nu_1)-x^\star(\nu_2),C^T(\nu_2-\nu_1)\rangle\\
&=\langle
x^\star(\nu_1)-x^\star(\nu_2),H^{-1}C^T(\nu_2-\nu_1)\rangle_H\\
&\leq \|
x^\star(\nu_1)-x^\star(\nu_2)\|_H\|H^{-1}C^T(\nu_2-\nu_1)\|_H\\
&\leq \|H^{-1}C^T(\nu_2-\nu_1)\|_H^2\\
&= (\nu_2-\nu_1)^TCH^{-1}C^T(\nu_2-\nu_1)\\
&= \|\nu_2-\nu_1\|_{CH^{-1}C^T}^2
\end{align*}
where the first inequality is due to Cauchy-Schwarz inequality and the
second comes from Lemma~\ref{lem:helpResult}.
Applying Corollary~\ref{cor:quadUBequivConcave} gives the result.
\end{pf}

Next, we show that if $f$ is a strongly convex quadratic function
and $h$ satisfies certain conditions, then
Theorem~\ref{thm:dualQuadBound} gives the best possible bound of the
form \eqref{eq:dualQuadBound}.
\begin{prop}
Assume that $f(x) = \tfrac{1}{2}x^THx+\zeta^Tx$ with
$H\in\mathbb{S}_{++}^n$ and $\zeta\in\mathbb{R}^n$ and that there exists a set
$\mathcal{X}\subseteq\mathbb{R}^n$ with non-empty interior on which
$h$ is linear, i.e. $h(x) =
\xi_{\mathcal{X}}^Tx+\theta_\mathcal{X}$ for all $x\in\mathcal{X}$. Further,
assume that there exists $\widetilde{\nu}$ such that
$x^\star(\widetilde{\nu})\in{\rm{int}}(\mathcal{X})$.
Then for any matrix $\mathbf{L}\not\succeq
CH^{-1}C^T$, there exist $\nu_1$ and $\nu_2$ such that
\eqref{eq:dualQuadBound} does not hold.
\label{prp:dualQuadBoundTight}
\end{prop}
\begin{pf}
Since
$x^\star(\widetilde{\nu})\in{\rm{int}}(\mathcal{X})$ we get 
for all $\nu_{\epsilon}\in\mathcal{B}_{\epsilon}^{m+p}(0)$, where the radius
$\epsilon$ is small enough, that
$x^\star(\widetilde{\nu})-H^{-1}C^T\nu_{\epsilon}\in\mathcal{X}$.
Introducing $x_{\epsilon}=-H^{-1}C^T\nu_{\epsilon}$, we get from the
optimality conditions to \eqref{eq:innerMinimizer}
(that specifies $x^{\star}(\nu)$) that
\begin{align*}
0&=Hx^\star(\widetilde{\nu})+\zeta+\xi_\mathcal{X}+C^T\widetilde{\nu}\\
&=H(x^\star(\widetilde{\nu})+x_{\epsilon})+\zeta+\xi_\mathcal{X}+C^T(\widetilde{\nu}+\nu_{\epsilon})\\
&=H(x^\star(\widetilde{\nu})+x_{\epsilon})+\zeta+h^{\prime}
(x^\star(\widetilde{\nu})+x_{\epsilon})+C^T(\widetilde{\nu}+\nu_{\epsilon})
\end{align*}
where $h^{\prime}(x^\star(\widetilde{\nu})\in\partial
h(x^\star(\widetilde{\nu})$ and
$x^\star(\widetilde{\nu})+x_{\epsilon}\in\mathcal{X}$ is used in
the last step. This implies that
$x^\star(\widetilde{\nu}+\nu_{\epsilon}) =
x^\star(\widetilde{\nu})+x_{\epsilon}$ and consequently that
$x^\star(\widetilde{\nu}+\nu_{\epsilon})\in\mathcal{X}$ for any
$\nu_{\epsilon}\in\mathcal{B}_{\epsilon}^{m+p}(0)$. Thus, for any
$\nu\in \widetilde{\nu}\oplus\mathcal{B}_{\epsilon}^{m+p}(0)$
we get
\begin{align*}
d(\nu) &= \min_{x}\tfrac{1}{2}x^THx+\zeta^Tx+h(x)+\nu^T(Cx-c)\\
&= \min_{x}\tfrac{1}{2}x^THx+\zeta^Tx+\xi_{\mathcal{X}}^Tx+\nu^T(Cx-c)\\
&=-\tfrac{1}{2}\nu^TCH^{-1}C^T\nu+\xi^T\nu+\theta
\end{align*}
where $\xi\in\mathbb{R}^n$ and $\theta\in\mathbb{R}$ collects the
linear and constant terms respectively. Since on the set
$\widetilde{\nu}\oplus\mathcal{B}_{\epsilon}^{m+p}(0)$,
$d$ is a quadratic with Hessian $CH^{-1}C^T$, it is straight-forward
to verify that \eqref{eq:dualQuadBound} holds
with equality for all
$\nu_1,\nu_2\in\widetilde{\nu}\oplus\mathcal{B}_{\epsilon}^{m+p}(0)$
if $\mathbf{L}=CH^{-1}C^T$. Thus, since
$\widetilde{\nu}\oplus\mathcal{B}_{\epsilon}^{m+p}(0)$ has non-empty
interior, we can for any matrix $\mathbf{L}\not\succeq
CH^{-1}C^T$ find
$\nu_1,\nu_2\in\widetilde{\nu}\oplus\mathcal{B}_{\epsilon}^{m+p}(0)$
such that 
\begin{equation*}
\|\nu_1-\nu_2\|_{CH^{-1}C^T}\geq \|\nu_1-\nu_2\|_{\mathbf{L}}.
\end{equation*}
This implies that for any $\mathbf{L}\not\succeq
CH^{-1}C^T$ there exist
$\nu_1,\nu_2\in\widetilde{\nu}\oplus\mathcal{B}_{\epsilon}^{m+p}(0)$
such that
\begin{align*}
d(\nu_1)&= d(\nu_2)+\langle \nabla
d(\nu_2),\nu_1-\nu_2\rangle-\tfrac{1}{2}\|\nu_1-\nu_2\|_{CH^{-1}C^T}\\
&\leq d(\nu_2)+\langle \nabla
d(\nu_2),\nu_1-\nu_2\rangle-\tfrac{1}{2}\|\nu_1-\nu_2\|_{\mathbf{L}}
\end{align*}
This concludes the proof.
\end{pf}
Proposition~\ref{prp:dualQuadBoundTight} shows that the bound in
Theorem~\ref{thm:dualQuadBound} is indeed the best obtainable bound of
the form \eqref{eq:dualQuadBound} if $f$ is a quadratic and $h$
specifies the stated assumptions. Examples of functions
that satisfy the assumptions on $h$ in
Proposition~\ref{prp:dualQuadBoundTight}  include linear functions, indicator
functions of closed convex constraint sets with non-empty interior, and the
1-norm. However, indicator functions for affine subspaces do not
satisfy the the assumptions of
Proposition~\ref{prp:dualQuadBoundTight} since their interiors are empty
(except for the trivial sub-space $\mathbb{R}^n$). In the following
proposition we will present a result that shows how
Theorem~\ref{eq:dualQuadBound} can be improved in that case.
\begin{prop}
Assume that $f(x)=\frac{1}{2}x^THx+\zeta^Tx$ with $H\in\mathbb{S}_{++}^n$
and $\zeta\in\mathbb{R}^n$, and
that $h=I_{Ax=b}$. Then \eqref{eq:dualQuadBound} holds for all
$\mathbf{L}\in\mathbb{S}_{+}^{m+p}$ such that 
$\mathbf{L}\succeq CH^{-1/2}(I-M)H^{-1/2}C^T$ where $M =
H^{-1/2}A^T(AH^{-1}A^T)^{-1}AH^{-1/2}$. Further, for any matrix
$\mathbf{L}\not\succeq CH^{-1/2}(I-M)H^{-1/2}C^T$ there exist $\nu_1,
\nu_2\in\mathbb{R}^{m+p}$ such that \eqref{eq:dualQuadBound} does not hold.
\label{prp:hAffine1}
\end{prop}
\begin{pf}
We have
\begin{align}
\nonumber d(\nu)&=-F^\star(-C^T\nu)-c^T\nu\\
\nonumber &=-\sup_x \left(-\nu^TCx-f(x)-h(x)\right)-c^T\nu\\
\label{eq:dualPrp} &= \inf_x \left(\nu^TCx+\tfrac{1}{2}x^THx+\zeta^Tx+I_{Ax=b}(x)\right)-c^T\nu
\end{align}
since $F=f+h$. The solution $x^\star(\nu)$ to the minimization problem satisfies the
following KKT-equations
\begin{equation}
\begin{bmatrix}
H & A^T\\
A& 0
\end{bmatrix}
\begin{bmatrix}
x^\star(\nu)\\
\lambda^\star(\nu)
\end{bmatrix} =
\begin{bmatrix}
-C^T\nu-\zeta\\
b
\end{bmatrix}
\label{eq:KKTeq}
\end{equation}
where $\lambda^\star(\nu)$ are dual variables corresponding to the equality
constraints. We have
\begin{align*}
x^\star(\nu) &= -H^{-1}(A^T\lambda^\star(\nu)+C^T\nu+\zeta).
\end{align*}
Inserting this into the second set of equations in \eqref{eq:KKTeq} gives
\begin{align*}
-AH^{-1}(A^T\lambda^\star(\nu)+C^T\nu+\zeta)=b.
\end{align*}
Since by assumption $A$ has full row rank and $H$ in positive
definite, $AH^{-1}A^T$ is invertible. Introducing the notation $H_{A}=
AH^{-1}A^T$, this implies that
\begin{align*}
\lambda^\star(\nu)=-H_A^{-1}(AH^{-1}(C^T\nu+\zeta)+b)
\end{align*}
which in turn implies that
\begin{align*}
x^\star(\nu) &= H^{-1}(A^TH_A^{-1}(AH^{-1}(C^T\nu+\zeta)+b)-C^T\nu-\zeta)\\
&=H^{-1}(A^TH_A^{-1}AH^{-1}-I)(C^T\nu+\zeta)+H^{-1}A^TH_A^{-1}b\\
&=-H^{-1/2}(I-M)H^{-1/2}(C^T\nu+\zeta)+H^{-1}A^TH_A^{-1}b.
\end{align*}
Insertion of this into \eqref{eq:dualPrp} gives after straight-forward
computations that
\begin{align*}
d(\nu)&=-\frac{1}{2}\nu^TCH^{-1/2}(I-M)H^{-1/2}C^T\nu+\xi^T\nu+\theta
\end{align*}
where $\xi\in\mathbb{R}^{m+p}$ and $\theta\in\mathbb{R}$ collect the
linear and constant terms respectively. This implies that $d$ is a
concave quadratic function with negative Hessian 
$CH^{-1}(I-M)H^{-1/2}C^T$. For concave quadratic functions, it is
straight-forward to verify that \eqref{eq:dualQuadBound} holds with
equality for all $\nu_1,\nu_2\in\mathbb{R}^{m+p}$ if $\mathbf{L}$ is chosen as the negative Hessian, i.e.
$\mathbf{L}=CH^{-1/2}(I-M)H^{-1/2}C^T$. This further implies, 
that for any $\mathbf{L}\not\succeq
CH^{-1/2}(I-M)H^{-1/2}C^T$  there exist
$\nu_1,\nu_2\in\mathbb{R}^{m+p}$ such that \eqref{eq:dualQuadBound}
does not hold. This concludes the proof.
\end{pf}
For the preceding result to hold, it is actually sufficient to assume that $f$
is strongly convex on the null-space of $A$ since this results in an
unique solution of $x^\star(\nu)$. The corresponding result is stated
in the following proposition.
\begin{prop}
Assume that $f(x)=\frac{1}{2}x^THx+\zeta^Tx$ with $H\in\mathbb{S}_{+}^n$
and $\zeta\in\mathbb{R}^n$, and
that $h=I_{Ax=b}$. Further assume $x^THx>0$ whenever $x\neq 0$ and
$Ax=0$, i.e. that $H$ is positive definite on the
null-space of $A$. Then \eqref{eq:dualQuadBound} holds for all
$\mathbf{L}\in\mathbb{S}_{+}^{m+p}$ such that 
$\mathbf{L}\succeq CK_{11}C^T$ where 
\begin{equation}
\begin{bmatrix}
K_{11}&K_{12}\\
K_{21}&K_{22}
\end{bmatrix}=
\begin{bmatrix}
H&A^T\\
A&0
\end{bmatrix}^{-1}.
\label{eq:KKTinv}
\end{equation}
Further, for any matrix
$\mathbf{L}\not\succeq CK_{11}C^T$ there exist $\nu_1,
\nu_2\in\mathbb{R}^{m+p}$ such that \eqref{eq:dualQuadBound} does not hold.
\label{prp:hAffine2}
\end{prop}
\begin{pf}
Since $H$ is positive definite on the null-space of $A$, the
KKT-matrix in \eqref{eq:KKTeq} is invertible and $\left[\begin{smallmatrix}
  K_{11}&K_{12}\\K_{21} & K_{22}\end{smallmatrix}\right]$ exists, see
\cite[p. 523]{Boyd2004}.
Equation~\eqref{eq:KKTinv} implies that the solution the the
KKT-system \eqref{eq:KKTeq} is given by
\begin{equation*}
\begin{bmatrix}
x^\star(\nu)\\
\lambda^\star(\nu)
\end{bmatrix}=
\begin{bmatrix}
K_{11}&K_{12}\\
K_{21}&K_{22}
\end{bmatrix}
\begin{bmatrix}
-C^T\nu-\zeta\\
b
\end{bmatrix}.
\end{equation*}
That is, $x^\star(\nu)=-K_{11}(C^T\nu+\zeta)+K_{12}b$. Inserting this into
\eqref{eq:dualPrp} gives
\begin{align*}
d(\nu)&=-\frac{1}{2}\nu^TC(2K_{11}-K_{11}HK_{11})C^T\nu+\xi^T\nu+\theta\\
&=-\frac{1}{2}\nu^TCK_{11}C^T\nu+\xi^T\nu+\theta
\end{align*}
where again $\xi\in\mathbb{R}^{m+p}$ and $\theta\in\mathbb{R}$ collect
the linear and constant terms, and where $K_{11}HK_{11} = K_{11}$ is used in
the second equality. This identity follows from the upper left block
of $\left[\begin{smallmatrix}
  K_{11}&K_{12}\\K_{21} & K_{22}\end{smallmatrix}\right]\left[\begin{smallmatrix}
  H&A^T\\A & 0\end{smallmatrix}\right]\left[\begin{smallmatrix}
  K_{11}&K_{12}\\K_{21} & K_{22}\end{smallmatrix}\right]=\left[\begin{smallmatrix}
  K_{11}&K_{12}\\K_{21} & K_{22}\end{smallmatrix}\right]$ and using
$K_{11}^TA=K_{11}A=AK_{11}=0$, where $AK_{11}=0$ follows from the lower left block of $\left[\begin{smallmatrix}
  H&A^T\\A & 0\end{smallmatrix}\right]\left[\begin{smallmatrix}
  K_{11}&K_{12}\\K_{21} & K_{22}\end{smallmatrix}\right]=\left[\begin{smallmatrix}
  I&0\\0 & I\end{smallmatrix}\right]$. This implies that $d$ is a concave and
quadratic function with negative Hessian $CK_{11}C^T$, which implies that
\eqref{eq:dualQuadBound} holds with equality for any
$\nu_1,\nu_2\in\mathbb{R}^{m+p}$ if
$\mathbf{L}=CK_{11}C^T$. This further implies, 
that for any $\mathbf{L}\not\succeq
CK_{11}C^T$  there exist
$\nu_1,\nu_2\in\mathbb{R}^{m+p}$ such that \eqref{eq:dualQuadBound}
does not hold. This concludes the proof.
\end{pf}

\begin{rem}
In the model predictive control context, the preceding result
implies that the quadratic cost matrix associated with inputs should be
positive definite, while the quadratic cost matrix associated with the states
need only be positive semi-definite.
\end{rem}

\section{Fast dual gradient methods}

In this section, we will describe generalized fast gradient methods
and show how they 
can be applied to solve the dual problem \eqref{eq:dualProb}.
Generalized fast gradient
methods can be applied to solve problems of the form
\begin{equation}
{\hbox{minimize }} \ell(x)+\psi(x)
\label{eq:minfP}
\end{equation}
where $x\in\mathbb{R}^n$,
$\psi~:~\mathbb{R}^n\to\mathbb{R}\cup\{\infty\}$ is proper, closed and
convex, $\ell~:~\mathbb{R}^n\to\mathbb{R}$ is convex, differentiable, and
satisfies 
\begin{equation}
\ell(x_1)\leq \ell(x_2)+\langle \nabla
\ell(x_2),x_1-x_2\rangle+\tfrac{1}{2}\|x_1-x_2\|_{\mathbf{L}}^2
\label{eq:quadUpperBound}
\end{equation}
for all $x_1,x_2\in\mathbb{R}^n$ and some
$\mathbf{L}\in\mathbb{S}_{++}^n$.
Before we state the algorithm, we define the generalized prox operator
\begin{equation}
{\rm{prox}}_{\psi}^{\mathbf{L}}(x) := \arg\min_y
\left\{\psi(y)+\tfrac{1}{2}\|y-x\|_{\mathbf{L}}^2\right\}
\label{eq:proxOpDef}
\end{equation}
and note that
\begin{align}
\nonumber &{\rm{prox}}_{\psi}^{\mathbf{L}}\left(x- \mathbf{L}^{-1}\nabla \ell(x)\right) = \\
\nonumber &=\arg\min_y
\left\{\tfrac{1}{2}\|y-x+\mathbf{L}^{-1}\nabla \ell(x)\|_{\mathbf{L}}^2+\psi(y)\right\}\\
&=\arg\min_y \left\{\ell(x)+\langle \nabla \ell(x),y-x\rangle+\tfrac{1}{2}\|y-x\|_{\mathbf{L}}^2+\psi(y)\right\}.
\label{eq:proxEquiv}
\end{align}
The generalized fast gradient method is stated below.

\medskip
\hrule
\smallskip
\begin{alg} ~\\\indent{\bf{Generalized fast gradient method}}
\smallskip\hrule\smallskip
\noindent Set: $y^1= x^0\in\mathbb{R}^n, t^1=1$\\
{\bf{For}} $k\geq 1$
\begin{itemize}
\item[] $x^{k} =
  {\rm{prox}}_{\psi}^{\mathbf{L}}\left(y^k-\mathbf{L}^{-1}\nabla \ell(y^k)\right)$
\item[] $t^{k+1} = \frac{1+\sqrt{1+4(t^k)^2}}{2}$
\item[] $y^{k+1} = x^k+\left(\frac{t^k-1}{t^{k+1}}\right)(x^k-x^{k-1})$
\end{itemize}
\smallskip\hrule
\label{alg:GFGM}
\end{alg}
\medskip

The standard fast gradient method as presented in
\cite{BecTab_FISTA:2009} is obtained by setting $\mathbf{L} = LI$ in
Algorithm~\ref{alg:GFGM}, 
where $L$ is the Lipschitz constant to $\nabla \ell$.
The main step of the fast gradient method is to perform a
prox-step, i.e., to minimize \eqref{eq:proxEquiv} which can be seen as
an approximation of the function $\ell+\psi$. For the standard fast gradient
method, $\ell$ is approximated with a quadratic upper bound that has the
same curvature, described by $L$, in all directions. If this quadratic
upper bound is a bad approximation of the
function to be minimized, slow convergence rate properties are
expected. The generalization to allow for a matrix $\mathbf{L}$ in the
algorithm allows for quadratic upper bounds with different curvature
in different directions. This enables for quadratic upper bounds that
much better approximate the function $\ell$ and consequently gives
improved convergence rate properties.

The generalized fast gradient method has a convergence rate of (see
\cite{WangmengTIP2011}) 
\begin{equation}
  \ell_{\psi}(x^k)-\ell_{\psi}(x^\star)\leq \frac{2\|x^\star-x^0\|_{\mathbf{L}}^2}{(k+1)^2}
\label{eq:convRate}
\end{equation}
where $\ell_{\psi}:= \ell+\psi$. The convergence rate of the standard fast gradient
method as given in \cite{BecTab_FISTA:2009}, is obtained by setting
$\mathbf{L}=LI$ in \eqref{eq:convRate}.

The objective here is to apply the generalized fast gradient method to
solve the dual problem \eqref{eq:dualProb}. By introducing
$\widetilde{g}(\nu)=g^\star([0~I]\nu)$,
the dual problem \eqref{eq:dualProb} can be expressed $\max_{\nu}
d(\nu)-\widetilde{g}(\nu)$, where
$d$ is defined in \eqref{eq:dFcnDef}.
As shown in Theorem~\ref{thm:dualQuadBound}, the function $-d$ satisfies the
properties required to apply generalized fast gradient methods. Namely
that \eqref{eq:quadUpperBound} holds for
any $\mathbf{L}\in\mathbb{S}_{+}^{m+p}$ such that $\mathbf{L}\succeq CH^{-1}C^T$.
Further, since $g$ is a
closed, proper, and convex function so is $g^\star$, see \cite[Theorem
12.2]{Rockafellar}, and by \cite[Theorem
5.7]{Rockafellar} so is $\widetilde{g}$. This implies
that generalized fast gradient methods, i.e. Algorithm~\ref{alg:GFGM},
can be used to solve the dual problem \eqref{eq:dualProb}.
We set $-d = \ell$ and $\widetilde{g}=\psi$, and restrict
$\mathbf{L}={\rm{blkdiag}}(\mathbf{L}_{\lambda},\mathbf{L}_{\mu})$ to get the following
algorithm.

\medskip
\hrule
\smallskip
\begin{alg}~\\\noindent {\bf{Generalized fast dual gradient method}}
\smallskip\hrule\smallskip
\noindent Set: $z^1= \lambda^0\in\mathbb{R}^m,
v^1=\mu^0\in\mathbb{R}^p, t^1=1$\\
{\bf{For}} $k\geq 1$
\begin{itemize}
\item[] $y^k = \arg\min_{x}\left\{f(x)+h(x)+(z^k)^TAx+(v^k)^TBx\right\}$
\item[] $\lambda^{k} = z^k+\mathbf{L}_{\lambda}^{-1}(Ay^k-b)$
\item[] $\mu^{k} =
  {\rm{prox}}_{g^\star}^{\mathbf{L}_{\mu}}(v^k+\mathbf{L}_{\mu}^{-1}By^k)$
\item[] $t^{k+1} = \frac{1+\sqrt{1+4(t^k)^2}}{2}$
\item[] $z^{k+1} = \lambda^k+\left(\frac{t^k-1}{t^{k+1}}\right)(\lambda^k-\lambda^{k-1})$
\item[] $v^{k+1} =
  \mu^k+\left(\frac{t^k-1}{t^{k+1}}\right)(\mu^k-\mu^{k-1})$
\end{itemize}
\smallskip\hrule
\label{alg:GFDGM}
\end{alg}
\medskip

where $y^k$ is the primal variable at iteration $k$ that is used to
help compute the gradient $\nabla
d(\nu^k)$ where $\nu^k = (z^k,v^k)$. To arrive
at the $\lambda^k$ and $\mu^k$ iterations, we let $\xi^k =
(\lambda^k,\mu^k)$,
and note that
\begin{align}
\label{eq:proxDecomposition} \xi^k&={\rm{prox}}_{\widetilde{g}}^{\mathbf{L}}\left(\nu^k+\mathbf{L}^{-1}\nabla d(\nu^k)\right)\\
\nonumber &=\arg\min_{\nu} \left\{\tfrac{1}{2}\|\nu-\nu^k-\mathbf{L}^{-1}\nabla
d(\nu^k)\|_{\mathbf{L}}^2+g^\star([0
~I]\nu)\right\}\\
\nonumber &= \left[\begin{array}{l}
\arg\min_{z} \left\{\tfrac{1}{2}\|z-z^k-\mathbf{L}_{\lambda}^{-1}\nabla_{z}
d(\nu^k)\|_{\mathbf{L}_{\lambda}}^2\right\}\\
\arg\min_{v} \big\{\tfrac{1}{2}\|v-v^k-\mathbf{L}_{\mu}^{-1}\nabla_{v}
d(\nu^k)\|_{\mathbf{L}_{\mu}}^2+g^\star(v)\big\}
\end{array}\right]\\
\nonumber &=\left[\begin{array}{l}
z^k+\mathbf{L}_{\lambda}^{-1}(Ay^k-b)\\
{\rm{prox}}_{g^\star}^{\mathbf{L}_{\mu}}(v^k+\mathbf{L}_{\mu}^{-1}By^k)
\end{array}\right].
\end{align}
In the following proposition we state the convergence rate properties
of Algorithm~\ref{alg:GFDGM}.
\begin{prop}
Suppose that Assumption~\ref{ass:probAss} holds. If $\mathbf{L} =
{\rm{blkdiag}}(\mathbf{L}_\lambda,\mathbf{L}_\mu)\in\mathbb{S}_{++}^{m+p}$
is chosen such that $\mathbf{L}\succeq CH^{-1}C^T$. Then Algorithm~\ref{alg:GFDGM}
converges with the rate
\begin{equation}
D(\nu^\star)-D(\nu^k)
\leq \frac{2 \left\|
\nu^\star - \nu^0\right\|_{\mathbf{L}}^2}{(k+1)^2}, \forall k\geq 1
\label{eq:convRateD}
\end{equation}
where $D=d-\widetilde{g}$ and $k$ is the iteration number.
\label{prp:convRate}
\end{prop}
\begin{pf}
Algorithm~\ref{alg:GFDGM} is Algorithm~\ref{alg:GFGM} applied to solve
the dual problem \eqref{eq:dualProb}. The convergence rate of Algorithm~\ref{alg:GFGM} is
given by
\eqref{eq:convRate} provided that the function to be minimized a sum
of one convex, differentiable function that satisfies \eqref{eq:quadUpperBound} and
one closed, proper, and convex function, see
\cite{WangmengTIP2011}. The discussion preceding the presentation of
Algorithm~\ref{alg:GFDGM} shows that the dual function to be optimized satisfies these
properties for any $\mathbf{L}\in\mathbb{S}_{++}^{m+p}$ that satisfies
$\mathbf{L}\succeq CH^{-1}C^T$. This concludes the proof.
\end{pf}
\begin{rem}
If $h=I_{Ax=b}$, the requirement on
$\mathbf{L}$ in Proposition~\ref{prp:convRate} changes according to
the results presented in 
Proposition~\ref{prp:hAffine1} and 
Proposition~\ref{prp:hAffine2}.
\end{rem}
\begin{rem}
By forming a specific running average of previous primal variables, it
is possible to prove a $O(1/k)$ convergence rate for the distance to the
primal variable optimum and a $O(1/k^2)$ convergence rate for
the worst case primal infeasibility, see \cite{PatrinosTAC2013}.
\end{rem}

For some choices of conjugate functions $g^\star$,
${\rm{prox}}_{g^\star}^{\mathbf{L}_{\mu}}(x)$ can be difficult 
to evaluate. For
standard prox operators (given by ${\rm{prox}}_{g^\star}^{I}(x))$,
Moreau decomposition \cite[Theorem 31.5]{Rockafellar} states that
\begin{equation*}
{\rm{prox}}_{g^\star}^{I}(x)+{\rm{prox}}_{g}^{I}(x)=x.
\end{equation*}
In the following proposition, we will generalize this result to hold for
the generalized prox-operator used here.
\begin{prop}
Assume that $g~:~\mathbb{R}^n\to\mathbb{R}$ is a proper, closed, and
convex function. Then
\begin{equation*}
{\rm{prox}}_{g^\star}^{\mathbf{L}}(x)+\mathbf{L}^{-1}{\rm{prox}}_{g}^{\mathbf{L}^{-1}}(\mathbf{L}x)=x
\end{equation*}
for every $x\in\mathbb{R}^n$ and any $\mathbf{L}\in\mathbb{S}_{++}^n$.
\label{prp:genMoreau}
\end{prop}
\begin{pf}
Optimality conditions for the prox operator \eqref{eq:proxOpDef} give that
$y={\rm{prox}}_{g^\star}^{\mathbf{L}}(x)$ if and only if
\begin{equation*}
0\in \partial g^\star(y)+\mathbf{L}(y-x)
\end{equation*}
Introducing $v=\mathbf{L}(x-y)$ gives $v\in\partial g^\star(y)$ which
is equivalent to $y\in\partial g(v)$ \cite[Corollary 23.5.1]{Rockafellar}. Since $y =
x-\mathbf{L}^{-1}v$ we have
\begin{equation*}
0\in \partial g(v)+(\mathbf{L}^{-1}v-x)
\end{equation*}
which is the optimality condition for $v =
{\rm{prox}}_{g}^{\mathbf{L}^{-1}}(\mathbf{L}x)$. This concludes the proof.
\end{pf}
\begin{rem}
If $g=I_{\mathcal{X}}$ where $I_{\mathcal{X}}$ is the indicator
function, then $g^\star$ is the support function. Evaluating the prox
operator \eqref{eq:proxOpDef}
with $g^\star$ being a support function is difficult. However, through
Proposition~\ref{prp:genMoreau}, this can be rewritten to only require
the a projection operation onto the set $\mathcal{X}$.
If $\mathcal{X}$ is a box constraint and $\mathbf{L}$ is diagonal, then
the projection becomes a max-operation and hence very cheap to implement.
\end{rem}
\begin{rem}
We are not restricted to have one
auxiliary term $g$ only. We can have any number of auxiliary terms
$g_i$ that all decompose according to the computations in
\eqref{eq:proxDecomposition}, i.e.,
we get one prox-operation in the algorithm for every auxiliary term $g_i$.
\end{rem}

\section{Choosing the $\mathbf{L}$-matrix}

\label{sec:Lmatrix}

From Theorem~\ref{thm:dualQuadBound} and
Proposition~\ref{prp:hAffine2}, we get that the $\mathbf{L}$-matrix
used in the quadratic lower bound in the algorithm should 
satisfy $\mathbf{L}\succeq CPC^T$, where $P=H^{-1}$ or $P=K_{11}$
depending on if the assumptions in
Theorem~\ref{thm:dualQuadBound} or
Proposition~\ref{prp:hAffine2} are satisfied. 
To get as fast convergence as possible, the approximation of
the function $d$ used in the algorithm should as accurately as
possible resemble the function $d$ itself. In view of Theorem~\ref{thm:dualQuadBound} and
Proposition~\ref{prp:hAffine2}, we want
$\mathbf{L}$ to be a close as
possible to $CPC^T$. Letting $\mathbf{L}=(D^{T}D)^{-1}$, we propose
to achieve this by minimizing the condition number of $DCPC^TD^T$,
subject to $I\succeq DCPC^TD^T$. If there are no structural
constraints on $\mathbf{L}$ and if $CPC^T$ has full rank, then
minimizing the condition number of $DCPC^TD^T$ gives
$\mathbf{L}=(D^TD)^{-1}=CPC^T$. However, this situation is quite
uncommon. First, we often have structural constraints on $\mathbf{L}$
that need to be taken into 
account. The most common such structural constraint is diagonal $\mathbf{L}$,
since for separable $g$, the complexity of computing
${\rm{prox}}_{g}^{\mathbf{L}}(x)$ is not increased compared to using
$\mathbf{L}=LI$. Sometimes, block-diagonal $\mathbf{L}$ can be used,
or in rare cases, full matrices $\mathbf{L}$. All these structural
constraints - diagonal, block-diagonal, and full - can be represented as
follows: let
$\mathbb{L}$ be a set of pairs $(i,j)$ for which $\mathbf{L}_{ij}$ may be
non-zero, then
\begin{align*}
\mathcal{L}=\{\mathbf{L}\in\mathbb{S}_{++}^{m+p}~|~&\mathbf{L}=(D^TD)^{-1},\\
& 
D\in\mathbb{R}^{(m+p)\times (m+p)} {\hbox{ invertible}},\\
&\mathbf{L}_{ij}=[\mathbf{L}^{-1}]_{ij}=[D^{-1}]_{ij}=0 \hbox{ if } (i,j)\notin\mathbb{L}\}.
\end{align*}
For instance, letting
$\mathbb{L}=\{(1,1),(2,2)\ldots,(m+p,m+p)\}$ restricts
$\mathbf{L}\in\mathcal{L}$ to be diagonal. A second issue that
hinders the choice of $\mathbf{L}=CPC^T$, is that $\mathbf{L}$ is
restricted to be positive definite, while $CPC^T$ is positive
definite only if $C$ has full row rank and if $P$ is positive definite.
When $CPC^T$ is not positive definite, we 
instead propose to minimize the ratio between the largest and smallest
non-zero eigenvalues (since the eigenvalues that are zero
cannot be changed). Letting
$\lambda_{1}(DCPC^TD^T)$ be the largest non-zero eigenvalue of $DCPC^TD^T$ and 
$\lambda_{r}(DCPC^TD^T)$ be the smallest non-zero eigenvalue of
$DCPC^TD^T$ (where if $r=m+p$ all eigenvalues are non-zero), the
proposed optimization problems can be written as
\begin{equation}
D=\arg\min_{(D^TD)^{-1}\in\mathcal{L}} \frac{\lambda_{1}(DCPC^TD^T)}{\lambda_{r}(DCPC^TD^T)}.
\label{eq:minCondNbr}
\end{equation}
Next we will show how to solve \eqref{eq:minCondNbr} in the following
three cases, which include all problem instances we will encounter:
\begin{enumerate}[(C1)]
\item $CPC^T\in\mathbb{S}_{++}^{(m+p)}$ 
\item $QC^TCQ^T\in\mathbb{S}_{++}^{q}$ where $P=Q^TQ$ and
  $Q\in\mathbb{R}^{q\times n}$
\item ${\rm{rank}}(QC^TCQ^T)={\rm{rank}}(CPC^T)<\min(m+p,q)$
\end{enumerate}
Before we present how to compute the optimal preconditioner in each of
the three cases, we state the following lemma.

\begin{lem}
For any matrix $A\in\mathbb{R}^{m\times n}$,
the non-zero eigenvalues of $A^TA$ equals the non-zero eigenvalues of
$AA^T$.
\label{lem:eigenValSymm}
\end{lem}
\begin{pf}
Without loss of generality, we assume that $m\leq n$ and that
${\rm{rank}}(A)=q\leq m$.
Let $A=U\Sigma V^T$, be the singular value decomposition of $A$,
where $U\in\mathbb{R}^{m\times m}$ and 
$V\in\mathbb{R}^{n\times n}$ are orthonormal, and 
\begin{equation*}
\Sigma = \begin{bmatrix}
\begin{bmatrix}
s_1 && & \\
&\ddots& &\\
&&s_q\\
&&& 0
\end{bmatrix}&0
\end{bmatrix}\in\mathbb{R}^{m\times n}.
\end{equation*}
This implies that $AA^T=U\Sigma
V^TV\Sigma^T U^T=U(\Sigma\Sigma^T)U^T$, or equivalently that
$(AA^T)U=U(\Sigma\Sigma^T)$, and that $A^TA=V\Sigma^T U^TU\Sigma
V^T=V\Sigma^T\Sigma V^T$, or equivalently that
$(A^TA)V=V(\Sigma^T\Sigma)$. That is, the eigenvalues to $AA^T$ are
given by the diagonal entries of $\Sigma\Sigma^T$ and the eigenvalues
to $A^TA$ are given by the diagonal entries of $\Sigma^T\Sigma$, i.e.
the non-zero eigenvalues of $AA^T$ and $A^TA$ coincide. This concludes
the proof.
\end{pf}

\subsection{Case 1}

\label{sec:C1}

We consider Case 1, i.e. C1. This is the case considered in
Theorem~\ref{thm:dualQuadBound} with $P=H^{-1}$ and an additional
rank assumption on $C$. 
\begin{prop}
Assume that $CPC^T\in\mathbb{S}_{++}^{(m+p)}$.
Then a matrix $D$ with $(D^TD)^{-1}\in\mathcal{L}$ that
minimizes the
ratio \eqref{eq:minCondNbr} can be computed by solving the
semi-definite program
\begin{equation}
\begin{tabular}{ll}
minimize & $t$\\
subject to & $tCPC^T\succeq \mathbf{L}$\\
& $CPC^T\preceq \mathbf{L}$\\
& $\mathbf{L}\in\mathcal{L}$
\end{tabular}
\label{eq:lmiC1}
\end{equation}
where $\mathbf{L}=(D^TD)^{-1}$. Further, $\mathbf{L}\succeq CPC^T$.
\end{prop}
\begin{pf}
Since $CPC^T$ has full rank, \eqref{eq:minCondNbr} is 
the condition number. Thus, according to \cite[Section 3.1]{BoydLMI},
\eqref{eq:lmiC1} can be solved in order to minimize
\eqref{eq:minCondNbr} . Further, the second
constraint implies that $\mathbf{L}\succeq CPC^T$.
\end{pf}

\subsection{Case 2}

\label{sec:C2}

Here, we show how to minimize \eqref{eq:minCondNbr} in the second
case, C2. This covers both Theorem~\ref{thm:dualQuadBound} (with
$P=H^{-1}$) and Proposition~\ref{prp:hAffine2} (with $P=K_{11}$) with
the additional assumption that $C$ is not wide and has full column rank.

\begin{prop}
Assume that $QC^TCQ^T\in\mathbb{S}_{++}^{q}$, where
$P\in\mathbb{S}_{+}^n$ is factorized as $P=Q^TQ$, where
$Q\in\mathbb{R}^{q\times n}$ has rank $q$. Then
a matrix $D$ with $(D^TD)^{-1}\in\mathcal{L}$ that
minimizes the
ratio \eqref{eq:minCondNbr}
can be computed by solving the
semi-definite program
\begin{equation}
\begin{tabular}{ll}
minimize &  $-t$\\
subject to & $QC^TMCQ^T\preceq I$\\
& $QC^TMCQ^T\succeq tI$\\
& $M\in\mathcal{L}$
\end{tabular}
\label{eq:lmiC2}
\end{equation}
where $M=(D^TD)$. Further $\mathbf{L}=(D^TD)^{-1}\succeq CPC^T$.
\end{prop}

\begin{pf}
Since $QC^TMCQ^T$ has full rank, we get from
Lemma~\ref{lem:eigenValSymm}, we get that minimizing the condition
number 
of $QC^TMCQ^T$ is equivalent to minimizing the ratio between
the largest and smallest non-zero eigenvalues of $DCPC^TD^T$, i.e.
equivalent to solving
\eqref{eq:minCondNbr}. From \cite[Section 3.1]{BoydLMI}, we get that
\eqref{eq:lmiC2} minimizes the condition number of
$QC^TMCQ^T$ i.e. it minimizes \eqref{eq:minCondNbr}. Further,
the first inequality implies through
Lemma~\ref{lem:eigenValSymm} that
$DCPC^TD^TQ\preceq I$, which is equivalent to that
$\mathbf{L}=(D^TD)^{-1}\succeq CPC^T$.
This concludes the proof.
\end{pf}

\subsection{Case 3}

\label{sec:C3}

Here, we consider Case C3, which covers the cases not included in
Cases C1 and C2. This covers, e.g. the situation in 
Proposition~\ref{prp:hAffine2} with additional assumptions on the rank
of $C$.
\begin{prop}
Assume that ${\rm{rank}}(QC^TCQ^T)=r$ with $r<\min(m+p,q)$
and that $P\in\mathbb{S}_{+}^n$ is factorized as $P=Q^TQ$, where
$Q\in\mathbb{R}^{q\times n}$ has rank $q$. Further, assume that
$\Phi\in\mathbb{R}^{q\times r}$ is an orthonormal basis for
$\mathcal{R}(QC^T)$. Then
a matrix $D$ with $(D^TD)^{-1}\in\mathcal{L}$ that
minimizes the
ratio \eqref{eq:minCondNbr} 
can be computed by solving the
semi-definite program
\begin{equation}
\begin{tabular}{ll}
minimize & $-t$\\
subject to & $QC^TMCQ^T\preceq I$\\
& $\Phi^TQC^TMCQ^T\Phi\succeq tI$\\
& $M\in\mathcal{L}$
\end{tabular}
\label{eq:lmiC3}
\end{equation}
where $M=(D^TD)$. Further, $\mathbf{L}=(D^TD)^{-1}\succeq CPC^T$.
\end{prop}
\begin{pf}
The first inequality in \eqref{eq:lmiC3} is by
Lemma~\ref{lem:eigenValSymm} equivalent to that $DCPC^TD^T\preceq I$,
i.e. $\lambda_{1}(DCPC^TD^T)\leq 1$.

To lower bound the smallest nonnegative eigenvalue, we need to search
in directions perpendicular to the null-space of
$QC^TMCQ^T$. We have
that
\begin{align*}
\mathcal{N}(QC^TMCQ^T)=\mathcal{N}(DCQ^T)=\mathcal{N}(CQ^T)\perp\mathcal{R}(QC^T)
\end{align*}
where the second equality holds since $D$ is assumed invertible. This
implies that we need to search in directions that span
$\mathcal{R}(QC^T)$.
Now, we have that $t\leq \lambda_{r}(QC^TMCQ^T)$ if and only if 
$0\leq x^T(QC^TMCQ^T-tI)x$ for all $x\in\mathcal{R}(QC^T)$.
This, in turn, is
equivalent to that 
\begin{equation}
\Phi^T(QC^TMCQ^T-tI)\Phi\in\mathbb{S}_{+}^r
\label{eq:smallestEigCond}
\end{equation}
where
$\Phi\in\mathbb{R}^{q\times r}$ is an orthonormal basis to
$\mathcal{R}(QC^T)$. Further, since $\Phi$ is an orthonormal basis,
i.e. $\Phi^T\Phi=I$,
\eqref{eq:smallestEigCond} is equivalent to
$\Phi^TQC^TMCQ^T\Phi\succeq tI$. This chain of
equivalences shows that
the second inequality in \eqref{eq:lmiC3} is equivalent to that
$\lambda_r(QC^TMCQ^T)\geq t$. Thus, by maximizing $t$ (or
equivalently minimizing $-t$) the ratio
\begin{equation*}
\lambda_1(QC^TMCQ^T)/\lambda_r(QC^TMCQ^T)\leq 1/t
\end{equation*}
is
minimized. From Lemma~\ref{lem:eigenValSymm} and the reasoning to the
proof of Case C2, we conclude that \eqref{eq:lmiC3} 
solves \eqref{eq:minCondNbr}. 

Further, the first inequality
implies through Lemma~\ref{lem:eigenValSymm} that $\mathbf{L}=(D^TD)^{-1}\succeq
CPC^T$. This concludes the proof.
\end{pf}
\begin{rem}
Note that if ${\rm{rank}}(QC^TCQ^T)=q$, then $\Phi=I$ is an orthonormal
basis to $\mathcal{R}(QC^T)$ and \eqref{eq:lmiC3} reduces to
\eqref{eq:lmiC2}. Thus, \eqref{eq:lmiC3} is a generalization of
\eqref{eq:lmiC2} to cover also the positive semi-definite case. A
similar generalization that reduces to \eqref{eq:lmiC1} in the 
positive definite case would rely on searching in directions perpendicular to
$\mathcal{N}(DCPC^TD^T)=\mathcal{N}(QC^TD^T)\perp\mathcal{R}(DCQ^T)$
to lower bound the smallest non-zero eigenvalue. This implies that the
search directions depend on the decision variables $D$, which makes
such a generalization more elaborate.
\end{rem}

\section{Model predictive control}

In this section, we pose some standard model predictive control
problems and show how they can be solved using the methods
presented in this paper. The resulting algorithms will have simple
arithmetic operations only which 
allows for easier implementation in embedded systems. We
also show how to choose the $\mathbf{L}$-matrix in each case.

\medskip
\begin{exmp} We consider MPC optimization problems of the form
\begin{equation*}
\begin{tabular}{lll}
minimize 
&\multicolumn{2}{l}{$\displaystyle\sum_{t=0}^{N-1}\tfrac{1}{2}\left(x_t^TQx_t+u_t^TRu_t\right)+\tfrac{1}{2}x_N^TQ_fx_N$}\\
subject to & $x_{t+1} = \Phi x_t+\Gamma u_t$, & $t=0,\ldots,N-1$\\
& $x_{\min}\leq x_t \leq x_{\max}$,& $t=0,\ldots,N$\\
& $u_{\min}\leq u_t \leq u_{\max}$,& $t=0,\ldots,N-1$\\
& $x_0 = \bar{x}$
\end{tabular}
\end{equation*}
where $\bar{x},x_t\in\mathbb{R}^{n_x}$, $u_t\in\mathbb{R}^{n_u}$,
$\Phi\in\mathbb{R}^{n_x\times n_x}$, $\Gamma\in\mathbb{R}^{n_x\times n_u}$
and $Q\in\mathbb{S}_{++}^{n_x}$, $R\in\mathbb{S}_{++}^{n_u}$, 
$Q_f\in\mathbb{S}_{++}^{n_x}$ are all diagonal.
Letting $y = (x_0,\ldots,x_N,u_0,\ldots,u_{N-1})$, this can be cast as
\begin{equation*}
\begin{tabular}{lll}
minimize 
&$\tfrac{1}{2}y^THy$\\
subject to & $Ay=b\bar{x}$\\
& $y_{\min}\leq y \leq y_{\max}$
\end{tabular}
\end{equation*}
where $H$, $A$, $b$, $y_{\min}$, and $y_{\max}$ are structured
according to $y$.
We choose $f(y) = \frac{1}{2}y^THy$, $g = 0$, and $h = I_{\mathcal{Y}}$ where
$I_{\mathcal{Y}}$ is the indicator function to
\begin{equation*}
\mathcal{Y} = \{y\in\mathbb{R}^{(N+1)n_x+Nn_u}~|~y_{\min}\leq y \leq y_{\max}\}.
\end{equation*}
This implicitly implies that we introduce dual variables $\lambda$ for the
equality constraints $Ay=b\bar{x}$. The algorithm becomes:
\begin{align}
\label{eq:algEqDual1} &y^k = \arg\min_y \left\{\tfrac{1}{2}y^THy+I_{\mathcal{Y}}(y)+z^TAx\right\}\\
\label{eq:algEqDual2} &\lambda^{k} = z^k+\mathbf{L}_{\lambda}^{-1}(Ay^k-b\bar{x})\\
\label{eq:algEqDual3} &t^{k+1} = \tfrac{1+\sqrt{1+4(t^k)^2}}{2}\\
\label{eq:algEqDual4} &z^{k+1} = \lambda^k+\left(\tfrac{t^k-1}{t^{k+1}}\right)(\lambda^k-\lambda^{k-1})
\end{align}
where the first step \eqref{eq:algEqDual1} can be implemented as
\begin{equation}
y^k = \max\left(\min\left(-H^{-1}A^Tz^k,y_{\max}\right),y_{\min}\right)
\label{eq:ykMaxMin}
\end{equation}
due to the structure of the problem.
The preceding section suggests that
$\mathbf{L}_{\lambda}=(D^TD)^{-1}\succeq AH^{-1}A^T$ should be
chosen such that $I\approx DAH^{-1}A^TD^T$. Since $A$ is sparse and
$H^{-1}$ is diagonal
due to the MPC problem formulation, $D$ can be chosen to get
equality in $I\approx DAH^{-1}A^TD^T$, i.e. we can choose
$\mathbf{L}_{\lambda}=(D^TD)^{-1}=AH^{-1}A^T$.
The algorithm requires
the computation of $\mathbf{L}_{\lambda}^{-1}z$, where $z=Ay^k-b\bar{x}$, in each
iteration. Since $\mathbf{L}_{\lambda}=AH^{-1}A$ is sparse, this can
efficiently be implemented by offline storing the sparse Cholesky
factorization $R^TR=S^T\mathbf{L}_{\lambda}S$, where $R$ is sparse and upper
triangular, and $S$ is a permutation matrix. The online computation of
$\mathbf{L}_{\lambda}^{-1}z$ then reduces to one forward and one
backward solve, which can be very efficiently implemented.

The algorithm in this example is a generalization of the algorithm in
\cite{RichterMMOR2013}, where the matrix $\mathbf{L}$ is chosen as
$\mathbf{L}=\|AH^{-1}A^T\|_2I$. In the numerical section we will see
that this generalization can significantly improve the convergence rate.
\label{ex:richterGen}
\end{exmp}

Next, we present an algorithm that works for arbitrary positive
definite cost matrices, and arbitrary linear constraints.

\medskip
\begin{exmp}
We consider MPC optimization problems of the form
\begin{equation*}
\begin{tabular}{lll}
minimize 
&\multicolumn{2}{l}{$\displaystyle\sum_{t=0}^{N-1}\frac{1}{2}\left(x_t^TQx_t+u_t^TRu_t\right)+\frac{1}{2}x_N^TQ_fx_N$}\\
subject to & $x_{t+1} = \Phi x_t+\Gamma u_t$, & $t=0,\ldots,N-1$\\
& $\underline{d}_x\leq B_xx_t \leq \bar{d}_x$,& $t=0,\ldots,N-1$\\
& $\underline{d}_u\leq B_uu_t\leq \bar{d}_u$,& $t=0,\ldots,N-1$\\
& $x_0 = \bar{x}, \underline{d}_N\leq B_Nx_N\leq \bar{d}_N$
\end{tabular}
\end{equation*}
where $\bar{x},x_t\in\mathbb{R}^{n_x}$, $u_t\in\mathbb{R}^{n_u}$,
$\Phi\in\mathbb{R}^{n_x\times n_x}$, $\Gamma\in\mathbb{R}^{n_x\times
  n_u}$, $B_x\in\mathbb{R}^{p_x\times n_x}$,
$B_u\in\mathbb{R}^{p_u\times n_x}$, $B_N\in\mathbb{R}^{p_N\times
  n_x}$, $\underline{d}_x,\bar{d}_x\in\mathbb{R}^{p_x}$,
$\underline{d}_u,\bar{d}_u\in\mathbb{R}^{p_u}$,
$\underline{d}_N,\bar{d}_N\in\mathbb{R}^{p_N}$, 
$Q\in\mathbb{S}_{++}^{n_x}$,$R\in\mathbb{S}_{++}^{n_u}$, and 
$Q_f\in\mathbb{S}_{++}^{n_x}$.
We let $y = (x_0,\ldots,x_{N},u_0,\ldots,u_{N-1})$ and
define $B={\rm{blkdiag}}(\bar{B}_x,B_N,\bar{B}_u)$ where $\bar{B}_x =
{\rm{blkdiag}}(B_x,\ldots,B_x)$ and $\bar{B}_u =
{\rm{blkdiag}}(B_u,\ldots,B_u)$. We also introduce
$\underline{d}=(\underline{d}_x,\ldots,\underline{d}_x,\underline{d}_N,\underline{d}_u,\ldots,\underline{d}_u)$
and $\bar{d}=(\bar{d}_x,\ldots,\bar{d}_x,\bar{d}_N,\bar{d}_u,\ldots,\bar{d}_u)$. This implies that
all inequality constraints are described by $\underline{d}\leq By\leq \bar{d}$.
Using this notation, the optimization problem can be rewritten as
\begin{equation*}
\begin{tabular}{lll}
minimize
&$\frac{1}{2}y^THy$\\
subject to & $Ay=b\bar{x}$\\
&$By=v$\\
&$\underline{d}\leq v\leq \bar{d}$
\end{tabular}
\end{equation*}
We let $f(y) = \frac{1}{2}y^THy$, $h=I_{Ay=b\bar{x}}$, and $g = I_{\mathcal{Y}}$ where
$\mathcal{Y} = \{y\in\mathbb{R}^{(N+1)n_x+Nn_u}~|~\underline{d}\leq y\leq \bar{d}\}$. Since
$h$ is the indicator function for the equality constraints
$Ay=b\bar{x}$, we do not need to introduce dual variables for those
constraints. However, we introduce dual variables $\mu$ for $By=v$.
Letting $H_A = AH^{-1}A^T$, the algorithm becomes
\begin{align}
\label{eq:algIneqDual1} & y^k = H^{-1}(A^TH_A^{-1}(AH^{-1}B^Tv^k+b\bar{x})-B^Tv^k)\\
\label{eq:algIneqDual2} & \mu^{k} = {\rm{prox}}_{g^\star}^{\mathbf{L}_{\mu}}(v^k+\mathbf{L}_{\mu}^{-1}By^k)\\
\label{eq:algIneqDual3} & t^{k+1} = \tfrac{1+\sqrt{1+4(t^k)^2}}{2}\\
\label{eq:algIneqDual4} & v^{k+1} =
  \mu^k+\left(\tfrac{t^k-1}{t^{k+1}}\right)(\mu^k-\mu^{k-1})
\end{align}
where the $y^k$ iterate follows from solving $\min_x
\big\{f(x)+I_{Ax=b\bar{x}}(x)+(v^k)^T Bx\big\}$. In an implementation, the
$y^k$-update can be implemented as in \eqref{eq:algIneqDual1}. Then, for
efficiency, the matrix multiplications should be computed offline and
stored for online use. Depending on the sparsity of $H$, $A$, and $B$,
it might be more efficient to use the KKT-system from which
\eqref{eq:algIneqDual1} is deduced, namely
\begin{equation*}
\begin{bmatrix}
H & A^T\\
A & 0
\end{bmatrix}
\begin{bmatrix}
y^k\\
\xi
\end{bmatrix}=
\begin{bmatrix}
-B^Tv^k\\
b\bar{x}
\end{bmatrix}.
\end{equation*}
Then, a sparse LDL-factorization of the KKT-matrix $\left[\begin{smallmatrix}
H & A^T\\
A & 0
\end{smallmatrix}\right]$ is computed offline for online use.
The online computational burden to compute the $y^k$-update then becomes one
forward and one backward solve. Whichever method that has the lower
number of flops should be chosen.

\begin{table*}
\centering
\caption{Comparison to other first-order methods, all implemented in
  MATLAB.}
\begin{tabular}{llrrrr}
 & & \multicolumn{2}{c}{~~exec time (ms)} &
\multicolumn{2}{c}{nbr iters}\\
Algorithm & Parameters & avg. & max & avg. & max\\
\hline
\eqref{eq:algEqDual1}-\eqref{eq:algEqDual4} & $\mathbf{L}_{\lambda} = AH^{-1}A^T$ & 2.3 & 12.1 & 21.7 & 102\\
\cite{RichterMMOR2013}& $\mathbf{L}_{\lambda} = \|AH^{-1}A^T\|_2I$ &
4713.9& 28411 & 50845 & 308210 \\
\eqref{eq:algIneqDual1}-\eqref{eq:algIneqDual4} & $\mathbf{L}_{\mu}$
comp. as in Sec.\ref{sec:Lmatrix} w. $P=K_{11}$ & 1.4 & 7.1 &
23.5 & 128\\
\eqref{eq:algIneqDual1}-\eqref{eq:algIneqDual4} & 
$\mathbf{L}_{\mu}$
comp. as in Sec.\ref{sec:Lmatrix} w. $P=H^{-1}$ & 1.2 & 5.8 &
20.0 & 105\\
\cite{PatrinosTAC2013} & $\mathbf{L}_{\mu} = \|BK_{11}B^T\|_2I$ & 98.5 & 673.0 &
1835.9 & 12686\\
\cite{PatrinosTAC2013} & $\mathbf{L}_{\mu} = \|BH^{-1}B^T\|_2I$ & 98.9 & 679.4 &
1850.1 & 12783\\
\cite{ODonoghueSplitting,JerezMHz} & $\rho=0.3$ &
193.9 & 920.6 & 3129.5 & 15037\\
\cite{ODonoghueSplitting,JerezMHz} & $\rho=3$ &
29.7 & 142.2 & 457.3 & 2179\\
\cite{ODonoghueSplitting,JerezMHz} & $\rho=30$ &
35.1 & 264.4 & 556.7 & 4194\\
\hline
\end{tabular}
\label{tab:numEvalMatlab}
\end{table*}

By restricting $\mathbf{L}_{\mu}$ to be diagonal, the second
step, i.e. \eqref{eq:algIneqDual2}, can be implemented as
\begin{equation*}
\mu^{k} = \min(v^k+\mathbf{L}_{\mu}^{-1}(By^k-\underline{d}),\max(v^k+\mathbf{L}_{\mu}^{-1}(By^k-\bar{d}),0)).
\end{equation*}
To get fast convergence, the diagonal $\mathbf{L}_{\mu}$ should be
computed as in Section~\ref{sec:Lmatrix}. Note that, in this
example, the matrix $P$ used in
Section~\ref{sec:Lmatrix} can be either $P=K_{11}$, where
$K_{11}$ is implicitly 
defined in \eqref{eq:KKTinv}, or $P=H^{-1}$. Since $K_{11}\preceq
H^{-1}$, the latter choice is expected to give a somewhat slower
convergence.

The splitting method used here is the same as the one used in
\cite{PatrinosTAC2013}. However, this is more general since we allow
for $\mathbf{L}_{\mu}$-matrices that are not a multiple of the
identity matrix. Also, the same splitting is used in \cite{ODonoghueSplitting,JerezMHz},
where ADMM (see \cite{BoydDistributed}) is used to solve the optimization problem.
\label{ex:patrinosGen}
\end{exmp}

\section{Numerical example}

The proposed algorithms are evaluated by applying them to the AFTI-16
aircraft model in \cite{Kapasouris,Bemporad97}. This problem is also a
tutorial example in the MPC toolbox in MATLAB. As in
\cite{Bemporad97} and the MPC toolbox tutorial, the continuous
time model from \cite{Kapasouris} 
is sampled using zero-order hold every 0.05 s. The system has four
states $x=(x_1,x_2,x_3,x_4)$, two outputs $y=(y_1,y_2)$, two
inputs $u=(u_1,u_2)$,
and obeys the following  
dynamics
\begin{align*}
x^+ &= \left[\begin{smallmatrix}
    0.999  & -3.008 &  -0.113 &  -1.608\\
   -0.000  &  0.986 &   0.048 &   0.000\\
    0.000  &  2.083 &   1.009 &  -0.000\\
    0.000  &  0.053 &   0.050 &   1.000
\end{smallmatrix}\right]x+\left[\begin{smallmatrix}
   -0.080  & -0.635\\
   -0.029  & -0.014\\
   -0.868  & -0.092\\
   -0.022  & -0.002
\end{smallmatrix}\right]u,\\
y &= \left[\begin{smallmatrix}
    0 & 1 & 0 & 0\\
    0 & 0 & 0 & 1
\end{smallmatrix}\right]x
\end{align*}
where $x^+$ denotes the state in the next time step. The dynamics,
input, and output matrices are denoted by $\Phi$, $\Gamma$, $C$
respectively, i.e. we have $x^+ = \Phi x+\Gamma u, y=Cx$.
The system is unstable, the magnitude of the largest eigenvalue 
of the dynamics matrix is 1.313. 
The outputs are the
attack and pitch angles, while the inputs are the elevator and
flaperon angles. The inputs are physically constrained to satisfy
$|u_i|\leq 25^\circ$, $i=1,2$. The outputs are soft constrained to
satisfy $-s_1-0.5\leq y_1\leq 0.5 +s_2$ and $-s_3-100\leq y_2\leq 100
+s_4$ respectively, where $s=(s_1,s_2,s_3,s_4)\geq 0$ are slack
variables. The cost in each time step is
\begin{equation*}
\ell(x,u,s)=\frac{1}{2}\big((x-x_r)^TQ(x-x_r)+u^TRu+s^TSs\big)
\end{equation*}
where $Q = C^TQ_yC+Q_x$, where $Q_y = 10^2I$ and $Q_x=
{\rm{diag}}(10^{-4},0,10^{-3},0)$, $x_r$ is such that $y_r=Cx_r$ where
$y_r$ is the output reference that can vary in each step, $R =
10^{-2}I$, and $S = 10^6I$. This 
gives condition number $10^{10}$ of the full cost matrix.
Further, the terminal cost is $Q$, and the control and prediction
horizon is $N=10$.
The numerical data in
Tables~\ref{tab:numEvalMatlab} and
\ref{tab:numEvalC} is obtained by following
a reference trajectory on
the output. The objective is to change the pitch angle from $0^\circ$
to $10^\circ$ and then back to $0^\circ$ while the angle of
attack satisfies the output constraints $-0.5^\circ \leq y_1 \leq
0.5^\circ$. The constraints on the angle of attack limits the rate on
how fast the pitch angle can be changed.

In Table~\ref{tab:numEvalMatlab}, the proposed algorithms
are evaluated by comparing them to other first 
order methods recently proposed in the literature for embedded model
predictive control, namely
\cite{RichterMMOR2013,PatrinosTAC2013,ODonoghueSplitting,JerezMHz}. In
Table~\ref{tab:numEvalC}, the execution time of a C implementation of
Algorithm~\ref{alg:GFDGM} is compared to the execution time of
FORCES, \cite{FORCES}, which is a C code generator for
MPC-problems, and to the commercial solver MOSEK. 

All algorithms in the comparison in Table~\ref{tab:numEvalMatlab} are
implemented in MATLAB, while the algorithms in
Table~\ref{tab:numEvalC} are implemented in C. Further, all
simulations are performed on a Linux machine using a single core
running at 2.9 GHz. To create an easily
transferable and fair
termination criterion, the optimal solution to each optimization
problem $y^\star$ is computed to high accuracy using an interior point solver.
Where applicable, the optimality condition is $\|y^k-y^\star\|_2/\|y^\star\|_2\leq
0.005$, where $y^k$ is the primal iterate in the algorithm. This
implies that a relative accuracy of 0.5$\%$ of the primal solution is required.

First, we discuss the results in Table~\ref{tab:numEvalMatlab}. The
algorithms in Example~\ref{ex:richterGen}, i.e.
\eqref{eq:algEqDual1}-\eqref{eq:algEqDual4}, and
Example~\ref{ex:patrinosGen}, i.e.
\eqref{eq:algIneqDual1}-\eqref{eq:algIneqDual4}, have been applied to this problem. Due to
the slack variables, \eqref{eq:ykMaxMin} cannot replace
\eqref{eq:algEqDual1} for the $y^k$ update. However, the
$y^k$ minimization is separable 
in the constraints and each of the projections can be solved by a
multi-parametric program with two regions. This is almost as
computationally inexpensive as the $y^k$ update in \eqref{eq:ykMaxMin}. Further,
we use $\mathbf{L}_{\lambda} = AH^{-1}A^T$. Algorithm \eqref{eq:algEqDual1}-\eqref{eq:algEqDual4}
is a generalization of
\cite{RichterMMOR2013}
that allows for general matrices $\mathbf{L}_{\lambda}$. The algorithm in
\cite{RichterMMOR2013} is obtained by setting $\mathbf{L}_{\lambda} =
\|AH^{-1}A^T\|_2I$. The
numerical evaluation in Table~\ref{tab:numEvalMatlab} reveals that this
generalization improves the execution time with more than three orders
of magnitude for this problem. 
The formulation in Example~\ref{ex:patrinosGen}, i.e.
\eqref{eq:algIneqDual1}-\eqref{eq:algIneqDual4}, directly covers this
MPC formulation with soft constraints. For this algorithm, we compute
$\mathbf{L}_{\mu}$ as in Section~\ref{sec:Lmatrix} using both
$P=K_{11}$ and $P=H^{-1}$.
The resulting algorithm is a generalization of the algorithm in
\cite{PatrinosTAC2013}. The algorithm in \cite{PatrinosTAC2013} is
given by setting $\mathbf{L}_{\mu} = \|BH^{-1}B^T\|_2I$ or
$\mathbf{L}_{\mu} = \|BK_{11}B^T\|_2I$ in the
iterations \eqref{eq:algIneqDual1}-\eqref{eq:algIneqDual4}.
Table~\ref{tab:numEvalMatlab} indicates that this
generalization improves the algorithm by one to two orders of
magnitude compared to \cite{PatrinosTAC2013}. Further,
\eqref{eq:algIneqDual1}-\eqref{eq:algIneqDual4} is based on
the same splitting as the method in 
\cite{ODonoghueSplitting,JerezMHz}. The difference is that here, the
problem is solved with a generalized dual gradient
method, while in \cite{ODonoghueSplitting,JerezMHz} it is solved using
ADMM. In ADMM, the
$\rho$-parameter need to be chosen. However, no exact guidelines are yet
known for this choice, and the performance of the algorithm often
relies heavily on this parameter. We compare our algorithm with ADMM
using the best $\rho$
that we found, $\rho=3$, and with one larger and one smaller $\rho$.
Table~\ref{tab:numEvalMatlab} reports that the execution time for our
method is one to two orders of
magnitude smaller (or more if the 
$\rho$-parameter in \cite{ODonoghueSplitting,JerezMHz} is chosen
suboptimally) than the algorithm proposed in
\cite{ODonoghueSplitting,JerezMHz}.

\begin{table}
\centering
\caption{Comparison to state-of-the-art solvers, all implemented in C.}
\begin{tabular}{llrr}
 & & \multicolumn{2}{c}{~~exec time (ms)} \\
Algorithm & Parameters & avg. & max \\
\hline
\eqref{eq:algEqDual1}-\eqref{eq:algEqDual4} & $\mathbf{L}_{\lambda} =
AH^{-1}A^T$ & 0.079 & 0.232 \\ 
\eqref{eq:algIneqDual1}-\eqref{eq:algIneqDual4} & $\mathbf{L}_{\mu}$
as in Sec.\ref{sec:Lmatrix} w. $P=H^{-1}$ & 0.061 & 0.196 \\ 
FORCES & - & 0.347& 0.592 \\
MOSEK & - &
4.9 & 5.4 \\
\hline
\end{tabular}
\label{tab:numEvalC}
\end{table}

In Table~\ref{tab:numEvalC}, we compare different solvers
implemented in C. For the algorithms
presented in this paper, we generate C code that take the reference
trajectory and the initial state as inputs. Compared to the
corresponding MATLAB implementations in Table~\ref{tab:numEvalMatlab},
the generated C code is more than 20 times faster.
These implementations are compared to FORCES and MOSEK.
FORCES, see \cite{FORCES}, is an optimized interior
point C code generator for MPC problems. The structure of the MPC problem
is exploited to significantly reduce the computational time when
solving the KKT-system in each iteration. The
comparison also includes MOSEK, which is a general commercial
QP-solver that does not have the advantage of generating code for this
specific problem beforehand. 
The numerical evaluation in
Table~\ref{tab:numEvalC} shows that our algorithms and FORCES, for both
of which C code is generated for this specific problem instance,
outperform the general purpose commercial C solver MOSEK with more
than one
order of magnitude. Further, Table~\ref{tab:numEvalC} reveals that our
two algorithms perform similarly and that they are at least two to three
times faster than FORCES.

\section{Conclusions}

We have proposed a generalization of dual fast gradient methods. This
generalization allows the algorithm to, in each iteration, minimize a quadratic upper
bound to the negative dual function with different curvature in
different directions. This is in contrast to the
standard fast dual gradient method where a quadratic upper bound to
the negative dual with the same curvature in all directions is
minimized in each iteration. This generalization is made possible by
the main contribution of this paper that characterizes the set of
matrices that can be used to describe a quadratic upper bound to the
negative dual function. The numerical evaluation on an ill-conditioned
aircraft problem reveals that the
proposed algorithms outperform several other MPC problem solvers recently
proposed in the literature.

\bibliography{references}

\end{document}